\documentclass[12pt]{article}

\usepackage [latin1]{inputenc}

\usepackage{amssymb}

\newcommand{\beqn}{\begin{eqnarray}}
\newcommand{\eeqn}{\end{eqnarray}}
\newcommand{\be}{\begin{equation}}
\newcommand{\ee}{\end{equation}}
\newcommand{\ba}{\begin{array}}
\newcommand{\ea}{\end{array}}
\newcommand{\R}{{\rm\bf R}}
\newcommand{\C}{{\rm\bf C}}
\newcommand{\pa}{\partial}

\newcommand{\re}{\ref}
\newcommand{\ci}{\cite}
\newcommand{\la}{\label}
\newcommand{\bfr}{\begin{flushright}}
\newcommand{\efr}{\end{flushright}}
\newcommand{\bfl}{\begin{flushleft}}
\newcommand{\efl}{\end{flushleft}}
\newcommand{\fr}{\frac}
\newcommand{\ov}{\overline}
\newcommand{\ti}{\tilde}
\newcommand{\st}{\stackrel}

\newcommand{\we}{\wedge}
\newcommand{\De}{\Delta}

\newcommand{\om}{\omega}
\newcommand{\Om}{\Omega}

\newcommand{\na}{\nabla}

\newcommand{\Br}{|\kern-.25em|\kern-.25em|}
\newcommand{\brr}{{|\kern-.15em|\kern-.15em|\kern-.15em}\,}
\newcommand{\ddd}{\st{.\kern-.07em.\kern-.07em.}}
\def\N{{\rm I\kern-.1567em N}}                              % Doppel-N
\def\R{{\rm I\kern-.1567em R}}                              % Doppel R
\def\C{{\rm C\kern-4.7pt                                    % Doppel C
\vrule height 7.7pt width 0.4pt depth -0.5pt \phantom {.}}}
\def\Z  {{\sf Z\kern-4.5pt Z}}                              % Doppel Z
\newtheorem{proposition}{\bf Proposition}[section]

\begin{document}

\renewcommand{\theequation}{\thesection.\arabic{equation}}
\newtheorem{theorem}{Theorem}[section]
\renewcommand{\thetheorem}{\arabic{section}.\arabic{theorem}}
\newtheorem{definition}[theorem]{Definition}
\newtheorem{deflem}[theorem]{Definition and Lemma}
\newtheorem{lemma}[theorem]{Lemma}
\newtheorem{example}[theorem]{Example}
\newtheorem{remark}[theorem]{Remark}
\newtheorem{remarks}[theorem]{Remarks}
\newtheorem{cor}[theorem]{Corollary}
\newtheorem{pro}[theorem]{Proposition}

\newcommand{\bd}{\begin{definition}}
\newcommand{\ed}{\end{definition}}
\newcommand{\bt}{\begin{theorem}}
\newcommand{\et}{\end{theorem}}
\newcommand{\bqt}{\begin{qtheorem}}
\newcommand{\eqt}{\end{qtheorem}}

\newcommand{\bp}{\begin{pro}}
\newcommand{\ep}{\end{pro}}

\newcommand{\bl}{\begin{lemma}}
\newcommand{\el}{\end{lemma}}
\newcommand{\bc}{\begin{cor}}
\newcommand{\ec}{\end{cor}}

\newcommand{\bex}{\begin{example}}
\newcommand{\eex}{\end{example}}
\newcommand{\bexs}{\begin{examples}}
\newcommand{\eexs}{\end{examples}}

\newcommand{\bexe}{\begin{exercice}}
\newcommand{\eexe}{\end{exercice}}

\newcommand{\br}{\begin{remark} }
\newcommand{\er}{\end{remark}}
\newcommand{\brs}{\begin{remarks}}
\newcommand{\ers}{\end{remarks}}

\newcommand{\pru}{{\bf Proof~~}}

%%%%%%%%%%%%%%%%%%%%%%%%%%%%%%%%%%%%%%%%%%%%%%%%%%
%%%%%%%%%%%%%%%%%%%%%%%%%%%%%%%%%%%%%%%%%%%%%%%%%%

\begin{center}
{\Large Valeriy M. Imaykin}
\\
{\small Moscow High School No 179}
\\
{\small Moscow, Russian Federation}
\\
{\small E-mail: ivm61@mail.ru}
\end{center}

%ÓÄÊ 517.955

\vskip10mm

\begin{center}
{\large\bf ON STABILITY OF SOLITONS AND THEIR ATTRACTION}

\vspace{0.5cm}

{\large\bf FOR A ROTATING CHARGE WITH FIXED}

\vspace{0.5cm}

{\large\bf  MASS CENTER IN THE MAXWELL FIELD}

\vspace{5mm}

\parbox{146mm}
{\noindent We consider the system of Maxwell equations and Lorentz torque equation which describes a motion of charge in electromagnetic field. Under certain symmetry conditions on charge distribution and on initial fields the mass center of the charge remains fixed and the charge rotates around it. The system admits stationary soliton-type solutions. We study the Lyapunov and the orbital stability of the solitons exploiting the energy conservation law.

We also show, by the angular momentum argument,  that there is no attraction to a soliton of finite angular momentum on the surface of states of the same angular momentum.

The bibliography: 23 refs.}

\end{center}

\section{Introduction}

\noindent
In this note we study some dynamical aspects of the Abraham model for a classical charge with spin coupled to the Maxwell field. The model is known since early 1900-ths, cf. \cite{Abr,Spo}, and was the subject of a number of papers, both at the physical and mathematical levels of rigor, see below the comment on previous works and results.

First let us write down the equations of motion. The Maxwell field consists of the electric field $E(x,t)$ and the magnetic field $B(x,t)$. The charge has the center of mass $q$ with the velocity $\dot q$. We assume that the mass distribution, $m\,\rho(x)$, and the charge distribution, $e\,\rho(x)$, are proportional to each other. Here $m$ is the total mass, $e$ is the total charge, and we use a system of units, where $m=1$ and $e=1$; $\rho(x)$ is a smooth radially symmetric smearing function of compact support, explicitly,
\begin{equation}\label{eq1}
\rho\in C_0^{\infty}(\R^3),\,\,\,\,\rho(x)= \rho_{rad}(|x|),\,\,\,\,\rho(x)=0\,\,\,{\rm for}\,\,\,|x|>R_{\rho}>0.%\eqno{(C)}
\end{equation}
The angular velocity of the charge is denoted by $\om(t)\in\R^3$. In detail, $\om(t)$ is the angular velocity ``in space'' (in
the terminology of \ci{AKN}) of the charge. Namely, let us fix a
``center'' point $O$ of the particle as a rigid body. Then the trajectory of each
particular point of the body is described by
$$
x(t)=q(t)+R(t)(x(0)-q(0)),
$$
where $q(t)$ is the position of $O$ at the time $t$, and $R(t)\in
SO(3)$. Respectively, the velocity reads
\be\la{dotx}
\dot x(t)=\dot q(t)+\dot R(t)(x(0)-q(0))=\dot q(t)+\dot
R(t)R^{-1}(t)(x(t)-q(t))=\dot q(t)+\om(t)\we(x(t)-q(t)),
\ee
where $\om(t)\in\R^3$ corresponds to the skew-symmetric matrix $\dot
R(t)R^{-1}(t)$ by the rule
\be\la{mv}
\dot
R(t)R^{-1}(t)={\cal J}\om(t):=\left(
\ba{ccc}
0 & -\om_3(t) & \om_2(t) \\
\om_3(t) & 0 & -\om_1(t) \\
-\om_2(t) & \om_1(t) & 0
\ea
\right).
\ee
We assume that $x$ and $q$ refer to a certain Euclidean coordinate
system in $\R^3$, and the vector product $\we$ is defined in this
system by standard formulas. The identification (\re{mv}) of a
 skew-symmetric matrix and the corresponding angular velocity
vector is true in any Euclidean coordinate system of the same
orientation as the initial one.

Then the Maxwell equations read, \cite{Spo}:
\begin{equation}\label{mls}
\dot E(x,t)=\na\we
B(x,t)-(\dot q(t)+\om(t)\we(x-q(t)))\rho(x-q(t)),\quad\dot B(x,t)= - \na\we E(x,t),
\end{equation}
where the current has a contribution also from the internal rotation; together with the constraints,
\begin{equation}\label{div}
\na\cdot E(x,t)=\rho (x-q(t)),\,\,\,\,\na\cdot B(x,t)= 0.
\end{equation}
The back reaction of the field onto the charge is given through the Lorentz force equation
\begin{equation}\label{lf}
\ddot q(t)=\int\, [E(x,t)+(\dot q(t)+\om(t)\we(x-q(t)))\we B(x,t)]\rho(x-q(t))\,dx
\end{equation}
and the Lorentz torque equation
\begin{equation}\label{ltq}
I\,\dot \om(t)= \int \, (x-q(t))\we[E(x,t)+(\dot q(t)+\om(t)\we(x-q(t)))\we B(x,t)]\rho(x-q(t))\,dx,
\end{equation}
with the moment of inertia
\begin{equation}\label{ib}
I=\frac23\int \, x^2\rho(x)d\,x.
\end{equation}

The important question is to obtain solutions having constant velocity and of the form
\begin{equation}\label{solx}
q(t)=q+vt,\,\om(t)\equiv\om,\,E(x,t)=E_{v,\om}(x-vt),\,B(x,t)=B_{v,\om}(x-vt).
\end{equation}
We will call them the ``soliton solutions'', in brief, the ``solitons''. If in (\ref{mls}), (\ref{lf}), (\ref{lt}) we set $\om=0$, by hand,  then for every $v\in\R^3$ there is a unique solution of the form (\ref{solx}). However, for the Abraham model including spin,
the equation
(\ref{lt}) can be satisfied only if either $\om\Vert v$ or $\om\bot v$ \cite{Spo}. This result is surprising at first sight,
but reflects the semirelativistic nature of the Abraham model. The velocity singles out a direction, which is then
taken by $\om$. Eventually one has to understand the domain of attraction of this soliton-like solutions. In this paper we restrict ourselves however to a simpler situation, where the charge remains at rest for all times, $q\equiv0$. This can be achieved by assuming the (anti-) symmetry conditions
\begin{equation}\label{odd-even}
E(-x)=-E(x),\quad B(-x)=B(x)
\end{equation}
for the initial fields.
Then this property would persist for all times:
\begin{equation}\label{symm}
E(-x,t)=-E(x,t),\,\,\,B(-x,t)=B(x,t).
\end{equation}

The Lorentz force equation is automatically satisfied, the Maxwell equations simplify to
\begin{equation}\label{M}
\dot E(x,t)=\na\we
B(x,t)-(\om(t)\we x)\rho(x),\quad\dot B(x,t)= - \na\we E(x,t),
\end{equation}
with the constraints
\begin{equation}\label{trans}
\na\cdot E(x,t)=\rho(x),\,\,\,\na\cdot B(x,t)=0,
\end{equation}
and the Lorentz torque equation simplifies to
\begin{equation}\label{lt}
I\dot\om(t)=\int x\we[E(x,t)+(\om(t)\we x)\we B(x,t)]\rho(x)\,dx.
\end{equation}
The system (\ref{symm})--(\ref{lt}) is the subject of study of the present paper. We will consider solutions of finite energy
\begin{equation}\label{en}
H(\om,E,B)=\frac{I\om^2}{2}+
\frac 12\int \,\Big(|E(x)|^2+|B(x)|^2\Big)\,dx
%K
<\infty.
%%%K
\end{equation}
The corresponding phase space will be equipped with a suitable topology below and the existence and uniqueness of finite energy solutions will be briefly explained. Note that the total momentum of the system
\begin{equation}\label{P}
P:=\int E(x,t)\we B(x,t)\,dx=0
\end{equation}
in view of the symmetry conditions (\ref{symm}).

The solitons for the system (\ref{symm})-(\ref{lt}) have the form
$$
E(x,t)=E_\om(x),\,\,\, B(x,t)=B_\om(x),\,\,\, \om(t)=\om=const\in\R^3.
$$
The solitons satisfy the stationary equations
\begin{equation}\label{symms}
E_\om(-x)=-E_\om(x),\,\,\,B_\om(-x)=B_\om(x),
\end{equation}
\begin{equation}\label{Ms}
\na\we B_\om(x)-(\om\we x)\rho(x)=0,\,\,\,\na\we E_\om(x)=0,
\end{equation}
\begin{equation}\label{transs}
\na\cdot E_\om(x)=\rho(x),\,\,\,\na\cdot B_\om(x)=0,
\end{equation}
\begin{equation}\label{lts}
\int x\we[E_\om(x)+(\om\we x)\we B_\om(x)]\rho\,dx=0.
\end{equation}

To write down the exact formulas of solitons let us specify the version of Fourier transform we use. The Fourier transform $F[f](k)=\hat f(k)$ of a function $f(x)$ reads
\begin{equation}\label{F}
\hat f(k):=(2\pi)^{-3/2}\int e^{-ikx}f(x)\,dx.
\end{equation}
Then by the Fourier transform
\begin{equation}\label{Fprop}
xf(x)\mapsto i\na_k\hat f(k),\,\,\,\na f(x)\mapsto ik\hat f(k).
\end{equation}
The Parseval equality holds,
\begin{equation}\label{Pars}
\int f(x)\overline{g(x)}\,dx = \int \hat f(k)\overline{\hat g(k)}\,dk.
\end{equation}

Note that the conditions (\re{eq1}) imply some special properties of the Fourier transform $\hat\rho(k)$ of $\rho(x)$. First,
\be\la{eq1f}
\hat\rho(k)\,\,\,{\rm is\,\,\,a\,\,\,real-valued\,\,\,radial\,\,\,function:}\,\,\, \ov{\hat\rho(k)}=\hat\rho(k),\,\,\,\hat\rho(k)=\rho_r(r),\,\,r:=|k|.
\ee
Second, as for the spatial decay of $\hat\rho(k)$,
\be\la{eq2f}
\hat\rho(k)\,\,\,{\rm is\,\,\,at\,\,\,least\,\,\,a\,\,\,fast\,\,\,decaying\,\,\,function\,\,\,of\,\,\,the\,\,\, Schwarz\,\,\,space}.
\ee

Now, in Fourier space, the soliton fields are expressed \cite{Spo,IKS04} by
\begin{equation}\label{solf}
\hat E_\om(k)=\frac{-ik\hat\rho(k)}{k^2},\,\,\,\hat B_\om(k)=-\frac{k\we(\om\we\na_k\hat\rho(k))}{k^2}.
\end{equation}

It follows from (\re{eq1f}) to (\ref{solf}) that
\be\la{solsl2}
\hat E_\om\in L^2, \hat B_\om\in L^2\,\,\,\, {\rm in}\,\,\, k-{\rm space,\,\,\, hence}\,\,\,\, E_\om\in L^2, B_\om\in L^2\,\,\,\, {\rm in}\,\,\, x-{\rm space}.
\ee

A {\it stability of solitons} of the system (\ref{M}), (\ref{lt}) together with the symmetry conditions (\ref{symm}) and the constraints (\ref{trans})  (= ``spinning charge of fixed mass center'') is the main result of our paper. It can be viewed as a partial result on the way to establishing the soliton asymptotics and scattering behavior of solutions to the system, see the following comment.

Let us comment on previous works. Note that, concerning the above system, there is a number of formal analytical results, but very few mathematically rigorous results on the qualitative behavior of the system's solutions.

The system (\ref{mls}) to (\ref{ltq}) is well known since Abraham's works \cite{Abr,Abr2}. The direct derivation of the conservation laws
from (\ref{ltq}) is presented by Kiessling in \cite{Kiess}. Soliton solutions to the system (\ref{mls}) to (\ref{ltq}) were computed first by Schwarzschild \cite{Sch}, see also the derivation in \cite{Spo}.

Some papers concern Lagrangian and Hamiltonian structure of the system. In \cite[Section 11]{Abr2} Abraham
computed the Lagrangian as integral of $-A_0\rho+ \vec A\cdot\vec j$ for standing rotating spherically
symmetric electron subject to external fields obeying very special symmetry conditions. In this case the Lagrangian
depends only on one variable $\om$, the angular velocity. However, derivation of the torque equation (\ref{ltq}) from the
variational {\it Hamilton's least action principle} remained an open question.
The main goal of the Nodvik's paper \cite{Nod} is a variational derivation
of the Lorentz-covariant dynamics for the relativistic rotating charged particle
in the Maxwell field, and the proof of the corresponding conservation laws. The system of Nodvick's equations is overdetermined, since they do not include rotational bare inertia. The situation was improved by Appel and Kiessling in \cite{AK}, where they develop the theory for the
relativistic rotating particle introducing a re-normalization limit. An invariant derivation of the non-relativistic Abraham equation (\ref{ltq}) from the Hamilton least action principle relying on the Poincare equations on the Lie group SO(3) was provided in \cite{IKS15}, in \cite{IKS17} it was shown that the Kiessling conserved quantities are the N\"oter invariants of the system's Lagrangian. Hamiltonian structure of the system was clarified in \cite{BIM14}.

The new interest for the rather old Abraham model is caused by the fact that a broad class of models of this type display soliton-type asymptotics and scattering behavior as it was discovered in recent years, see e.g. \cite{IKM04,IKV12,IKS11}, where a charged particle moves in Maxwell or scalar field without rotation. In particular, in \cite{IKM04} the orbital stability was established for solitons of the Maxwell-Lorentz system for a moving but non-rotating particle. The method is a thorough combination of the energy and the total momentum conservation, for the system written in Hamiltonian form. Lagrangian and Hamiltonian structure of the models play a significant role in these methods, so it was of a considerable interest and importance, to include the Abraham model with rotating charge into the class of Lagrangian and Hamiltonian systems. Nevertheless, in view of rather complicated structure of the system, the problem of establishing soliton-type asymptotics and scattering behavior for the Abraham model with rotating charge remains open.

Some progress in this direction was made in the paper \cite{IKS04}, where results on soliton-type asymptotics in local energy seminorms and also on scattering of solitons in global energy norms are obtained for the system (\ref{symm}) to (\ref{lt}). The crucial assumption of the paper is that the norm of $\rho$ in $L^2$ is sufficiently small  that means a {\it weak wave-particle interaction}.

The method of the present paper does not need this assumption. It exploits {\it energy conservation arguments}. %and angular momentum

Using these arguments we prove, first, that the zero soliton with $\om=0$ is Lyapunov stable and orbital stable. Second, we prove stability of an arbitrary soliton with respect to a special class of perturbations of initial data, namely, perturbations of uniformly compact support, see Theorem \ref{t1} below\footnote{These two results are published in \ci{I21}. We present them for completeness, and also correct the computational error that was made in \ci{I21}, which did not affect the truth of its main result.}.

Recently the result on stability for the system (\ref{symm}) to (\ref{lt}) was improved in \ci{KKMil}, where the stability of the solitons is proved under the condition $I_{\rm eff}\gg I$, where $I$ opposite to (\ref{ib}) is considered as an independent parameter and $I_{\rm eff}:= I+(2/3)\int(|\na\hat\rho(k)|)/(k^2)\,dk$. The condition means that $\rho$ {\it is sufficiently large}.

Finally, in \ci{KKmov} it is claimed the general result on the stability of the solitons for rotating and moving particle, i.e. for the system (\ref{mls}) to (\ref{ib}), without any additional assumptions. Now the problem on stability can be considered as solved.

However, our approach deserves attention because it relies on a simple technique for using the law of conservation of energy. Moreover, the condition of uniform boundness of the support of field perturbations is physically meaningful.

Note that both papers \ci{KKMil} and \ci{KKmov} exploit the Hamilton structure of corresponding systems.

Then it occurs the question of {\it global attraction} of an arbitrary solution to the set of all solitons.

Global attraction of an arbitrary solution to the set of all solitons for 2D Maxwell-Lorentz system, in sufficiently weak weighted Sobolev norms, was proved in \ci{KKrest2D}.

In the present paper we establish a partial negative result for both 3D and 2D cases: there is no attraction, in energy norm, to a soliton of finite angular momentum, on the surface of states of the same angular momentum.

%%%%%%%%%%%%%%%%%%%%%%%%%%%%%%%%%%%%%%%%%%%%%%
%%%%%%%%%%%%%%%%%%%%%%%%%%%%%%%%%%%%%%%%%%%%%%

\section{Existence of dynamics and stability of solitons for special perturbations %the main results
}

\setcounter{equation}{0}

Consider the space $L = (\R^3, L^2(\R^3;\R^3), L^2(\R^3;\R^3))$ with the norm
\be\la{Lnorm}
\Vert(\om,E,B)\Vert_L:=|\om|+\Vert E\Vert_{L^2(\R^3;\R^3)}+\Vert B\Vert_{L^2(\R^3;\R^3)}.
\ee
\begin{definition}
The {\it phase space for the system} (\ref{symm})-(\ref{lt}) is the nonlinear submanifold ${\cal M}$ of states $(\om,E,B)\in L$, where  $E, B$ satisfy (\ref{symm}) and (\ref{trans}). The {\it topology of} ${\cal M}$ is defined through the imbedding ${\cal M}\subset L$.
\end{definition}

Let us state a proposition on existence of dynamics for the system (\ref{symm})-(\ref{lt}) in the introduced phase space ${\cal M}$.

\begin{proposition}\label{p1}
Let the assumptions (\ref{eq1}) hold. Then

i) for any $(\om_0,E_0,B_0)\in{\cal M}$ the Cauchy problem for the system (\ref{symm})-(\ref{lt}) has a unique solution $(\om(t),E(x,t),B(x,t))\in C(\R;{\cal M})$ with initial conditions $\om(0)=\om_0, E(x,0)=E_0,B(x,0)=B_0$;

ii) for any $T\in\R$ the map $U_T: (\om_0,E_0,B_0)\mapsto(\om(T),E(\cdot,T),B(\cdot,T))$ is continuous in ${\cal M}$;

iii) the energy $H(t):=H(\om(t),E(\cdot,t),B(\cdot,t))$ is conserved along the solutions of the system:
\begin{equation}\label{econ}
H(t)\equiv H(0),\,\,\,\,t\in\R.
\end{equation}

\end{proposition}

The proof is similar to that of \cite[Appendix]{KS00}, with the charge density $\rho(x)$ and the current $j(x,t):=(\om(t)\we x)\rho(x)$.

\smallskip

Now let us specify common definitions of Lyapunov stability, orbital stability, and asymptotic stability for the particular case of soliton solutions to the system (\ref{symm})-(\ref{lt}).

\begin{definition}
1) A soliton solution $Y_\om=(\om,E_\om,B_\om)$ is called {\it Lyapunov stable} if $\forall \varepsilon>0$ $\exists \delta >0$ such than for any solution $Y(t)=(\om(t),E_(t),B(t))$ from the condition $\Vert Y(0)-Y_\om(0)\Vert_L<\delta$ it follows that $\Vert Y(t)-Y_\om(t)\Vert_L<\varepsilon$ $\forall t\in\R$.

2) A soliton solution $Y_\om=(\om,E_\om,B_\om)$ is called {\it orbital stable} if $\forall \varepsilon>0$ $\exists \delta >0$ such than for any solution $Y(t)=(\om(t),E_(t),B(t))$ with $\Vert Y(0)-Y_\om(0)\Vert_L<\delta$ the solution $Y(\cdot)$ remains in $\varepsilon$-neighbourhood of $Y_\om(\cdot)$ in the space $C(\R;L)$.

3) A soliton solution $Y_\om=(\om,E_\om,B_\om)$ is called {\it asymptotic stable} if it is Lyapunov stable and $\Vert Y(t)-Y_\om(t)\Vert_L\to0$ as $t\to\infty$.
\end{definition}

Note that in view of the energy conservation (\ref{econ}) one could not expect the asymptotic stability of the solitons in the phase space ${\cal M}$ of finite energy states. Our main result is then  the Lyapunov and orbital stability. In fact, for our case the properties of Lyapunov stability and orbital stability of a soliton coincide, since solitons do not depend on $t$ and thus are stationary solutions (fixed points of the dynamical system).

Further, the method of exploiting energy conservation arguments we use imposes some restrictions on the class of possible perturbations. Namely, we assume that the perturbations are of {\it uniformly compact supports}, see the exact formulation below. The main result on stability is the following theorem.

\begin{theorem}\label{t1} a) The zero soliton with $\om=0$ is Lyapunov stable and orbital stable.

b) Let us fix a non-zero $\om\in\R^3$. Consider solutions to the Cauchy problem for the system (\ref{symm})-(\ref{lt}) with initial data $\om+\Omega_0,\,\,E_{\om}(x)+e_0(x),\,\,B_{\om}(x)+b_0(x)$, where $(\Om_0,e_0,b_0)\in L$, $e_0(x)$ is odd, $b_0(x)$ is even, $\na\cdot e_0=0$, $\na\cdot b_0=0$.

The soliton $(\om,E_{\om}(x),B_{\om}(x))$ is Lyapunov stable and as well orbital stable with respect to perturbations of uniformly compact support

i.e. $\forall \varepsilon>0$ and $\forall R>0$ $\exists \delta >0$ such that for any solution $(\om(t),E(t),B(t))$ with the initial data $\om_0=\om+\Omega_0, E_0=E_\om+e_0, B_0=B_\om+b_0$ from the condition $\Vert(\Omega_0, e_0, b_0)\Vert_L<\delta$ it follows that $\Vert (\om(t)-\om, E(t)-E_\om, B(t)-B_\om)\Vert_L<\varepsilon$ $\forall t\in\R$ $\forall (\Omega_0,e_0(x),b_0(x))$ such that ${\rm supp}\,e_0\subset\{|x|\le R\}$ and ${\rm supp}\,b_0\subset\{|x|\le R\}$.
\end{theorem}

%%%%%%%%%%%%%%%%%%%%%%%%%%%%%%%%%%%
%%%%%%%%%%%%%%%%%%%%%%%%%%%%%%%%%%%

\section{Proof of Theorem \ref{t1}}

\setcounter{equation}{0}

\subsection{Equations for perturbations %and their solutions
}

Let us fix an arbitrary soliton $(\om,E_\om,B_\om)$.
To study its stability we put
\begin{equation}\label{eb}
\om(t)=\om+\Om(t),\,\,E(x,t)=E_{\om}(x)+e(x,t),\,\,B(x,t)=B_{\om}(x)+b(x,t)
\end{equation}
with
\begin{equation}\label{transeb}
\na\cdot e = 0,\,\, \na\cdot b = 0,\,\, e(-x,t) = -e(x,t),\,\, b(-x,t) = b(x,t).
\end{equation}

Insert (\ref{eb}) into (\ref{M}) and (\ref{lt}), take the stationary equations (\ref{Ms})--(\ref{lts}) into account and obtain the following system for the perturbations $e$, $b$, $\Om$:
\begin{equation}\label{eqseb}
\dot e = \na\we b-(\Om\we x)\rho,\,\,\, \dot b = -\na\we e,
\end{equation}
\begin{equation}\label{eqOm}
I\dot\Om = \int x\we[e+(\Om\we x)\we B_{\om} +(\om\we x)\we b+(\Om\we x)\we b]\rho\, dx.	
\end{equation}

The following remark is very important for our further analysis.

\begin{remark}\label{r1}
{\rm i) The Cauchy problem for the system (\re{transeb}) -- (\re{eqOm}) with initial data $(\Om_0, e_0, b_0)$ apriori has the solution}
$$
\Om(t)=\om(t)-\om,\,\,\,e(x,t)=E(x,t)-E_\om(x),\,\,\,b(x,t)=B(x,t)-B_\om(x),
$$
{\rm where $(\om,E,B)$ it the solution to the system (\ref{symm})-(\ref{lt}) with the initial data $(\om+\Om_0, E_\om+e_0, B_\om+b_0)$.}

{\rm ii) By the energy conservation (\ref{econ}), (\ref{solsl2}), and (\ref{eb}),}
\begin{equation}\label{ub}
(\Om(t), e(x,t), b(x,t))\,\,\,{\rm is\,\, bounded\,\, in}\,\,\,L\,\,\,{\rm uniformly\,\,in}\,\,\,t\in\R.
\end{equation}
\end{remark}

Let us study the system  (\ref{transeb}), (\ref{eqseb}), (\ref{eqOm}).

First we express the fields $(e,b)$ from the system
\begin{equation}\label{eqsebj}
\dot e = \na\we b-j(x,t),\,\,\, \dot b = -\na\we e
\end{equation}
with zero charge density and the prescribed current $j(x,t)$. In our case
\be\la{cur}
j(x,t)=(\Om(t)\we x)\rho(x),\,\,\,\,\hat j(k,t)=(\Om(t)\we i\na_k)\hat\rho(k).
\ee

%%%%%%%%%%%%%%%%%%%%%%%%%%%%%%%%%%%%%%%%%%%%%%%%%%%%%%%%%%%%%%%%%%%%%%%%%%%%%%%
%%%%%%%%%%%%%%%%%%%%%%%%%%%%%%%%%%%%%%%%%%%%%%%%%%%%%%%%%%%%%%%%%%%%%%%%%%%%%%%%

One has
\begin{equation}\label{U}
\left(\ba{c} e(x,t) \\ b(x,t) \ea\right) = U(t)\left(\ba{c} e_0(x) \\ b_0(x) \ea\right)-
\int_0^tU(t-s)\left(\ba{c} j(x,s) \\ 0 \ea\right)\,ds,
\end{equation}
where $U(t)$ is the group of the free Maxwell equation. Note that the group is isometric in the space $[L^2(\R^3;\R^3)]^2$  by the corresponding energy conservation law for free Maxwell equations, \cite{IKS11,KS00,IKS02,KS06}. Put
\begin{equation}\label{zero}
\left(\ba{c} e_{(0)}(x,t) \\ b_{(0)}(x,t) \ea\right):=U(t)\left(\ba{c} e_0(x) \\ b_0(x) \ea\right),\,\,\,\, \left(\ba{c} e_{(r)}(x,t) \\ b_{(r)}(x,t) \ea\right):= -\int_0^tU(t-s)\left(\ba{c} j(x,s) \\ 0 \ea\right)\,ds.
\end{equation}

By \cite[Appendix]{KS00} in the Fourier space we get
\begin{equation}\label{ek}
\hat e=\frac{d}{dt}\hat K_t\hat e_0+im\hat K_t\hat b_0-\int_0^t \frac{d}{dt}\hat K\vert_{t-s}\hat j(s)\,ds,
\end{equation}
\begin{equation}\label{bk}
\hat b=-im\hat K_t\hat e_0+\frac{d}{dt}\hat K_t\hat b_0+\int_0^t im\hat K_{t-s}\hat j(s)\,ds.
\end{equation}
Here
\be\la{kernels}
m:=k\we,\,\,\, \hat K_t(k):=\frac{\sin(|k|t)}{|k|},\,\,\,\frac{d}{dt}\hat K(t)=\cos(|k|t).
\ee

For $\Om$ we obtain the closed equation
\begin{equation}\label{eqOmeb}
I\dot\Om = \int x\we[e+(\Om\we x)\we B_{\om} +(\om\we x)\we b+(\Om\we x)\we b]\rho\, dx,	
\end{equation}
where $e,b$ are given by (\ref{U}).

For $\om=0$ the equation (\re{eqOmeb}) reads
\begin{equation}\label{eqOmeb0}
I\dot\Om = \int x\we[e+(\Om\we x)\we b]\rho\, dx.	
\end{equation}
The system (\re{eqseb}), (\re{eqOmeb0}) is of the same type as the initial system (\re{odd-even}) to (\re{lt}) and the corresponding energy
$$
{\cal H}(t):=\frac{I\Om(t)^2}{2}+\frac 12\int \,\Big(|e(x,t)|^2+|b(x,t)|^2\Big)\,dx
$$
is conserved. This implies the statement a) of Theorem \re{t1}.

Further, we rewrite the equation (\ref{eqOm}) as
\begin{equation}\label{eqOm2}
I\dot\Om = \int x\we[(\Om\we x)\we B_{\om}+e+(\om(t)\we x)\we b]\rho\, dx.	
\end{equation}
In (\ref{eqOm2}), we consider $\om(t)$ as the known function, the first component of the solution to the system (\ref{symm}) -- (\ref{lt}) with the initial data $(\om+\Om_0,E_\om+e_0,B_\om+b_0)$. In particular, $\om(t)$ is uniformly bounded in $t\in\R$ by the conservation of the energy (\ref{econ}).

Below we derive, from equation (\re{eqOm2}), the stability of zero solution $(\Om,e,b)=(0,0,0)$ for an arbitrary finite time interval. The zero solution exists for the zero initial data $(\Om_0,e_0,b_0)=(0,0,0)$. At first we obtain a more detailed structure of equation (\re{eqOm2}).

\subsection{The structure of the equation (\ref{eqOm2})}

We write the right-hand side of (\ref{eqOm2}) as $T_1+T_2+T_3$ with
$$
T_1:=\int x\we[(\Om\we x)\we B_{\om}]\rho\, dx,\,\,\,T_2:=\int (x\we e)\rho\, dx,\,\,\,T_3:=\int x\we[(\om(t)\we x)\we b]\rho\, dx.
$$

{\it Step 1.} For the first term we have
$$
T_1=\int (\Om\we x)(x\cdot B_{\om})\rho\, dx=K\we\Om,
$$
where
\begin{equation}\label{K}
K:=-\int x(x\cdot B_{\om})\rho\, dx
\end{equation}
is a constant vector in $\R^3$. Let us compute $K$ in detail. Since $\hat\rho$ is real-valued,
$$
K=-\int (x\cdot B_{\om})x\rho\, dx = -\int (i\na_k\cdot \hat B_{\om})\overline{i\na_k\hat\rho}\, dk = -\int \na_k\cdot \hat B_{\om}\cdot\overline{\na_k\hat\rho}\, dk = \int\hat B_{\om}\Delta_k\hat\rho\, dk.
$$
Recall that
$$
\hat B_{\om} = -\frac{k\we(\om\we\na_k\hat\rho)}{k^2} = - \frac{\om(k\cdot\na_k\hat\rho)-\na_k\hat\rho(k\om)}{k^2} = \frac{k\tilde\rho(k\om)}{k^2} - \om\tilde\rho.
$$
Here we denote $r=|k|$, $\hat\rho=\rho_r(r)$, then $\na_k\hat\rho=k\rho_r'(r)/r$ and we put $\tilde\rho(r):=\rho_r'(r)/r$. Then for $K$ we obtain
$$
K=K_1-K_2;\,\,\,\,K_1:=\int\frac{\tilde\rho}{k^2}(k\otimes k)\om\,\Delta_k\hat\rho\,dk;\,\,\,\, K_2:=\om\int\,\tilde\rho\,\Delta_k\hat\rho\,dk.
$$
Further, by the skew-symmetry property of the non-diagonal terms of the matrix $k\otimes k$ we have
$$
K_1=\int\frac{\tilde\rho}{k^2}\Delta_k\hat\rho\left(\ba{c}k_1^2\,\om_1 \\ k_2^2\,\om_2 \\ k_3^2\,\om_3\ea\right)dk.
$$
Each of the integrals
$$
\int\frac{\tilde\rho}{k^2}\Delta_k\hat\rho\, k_j^2\om_j\,dk\,\,\,\,\,{\rm equals}\,\,\,\,\,\frac13\int\frac{\tilde\rho}{k^2}\Delta_k\hat\rho k^2\,\om_j\,dk= \frac13\int\tilde\rho\Delta_k\hat\rho \,\om_j\,dk,
$$
hence,
$$
K_1=\frac13{\om}\int\tilde\rho\Delta_k\hat\rho\,dk.
$$
Finally
$$
K_1-K_2=-\frac23{\om}\int\tilde\rho\Delta_k\hat\rho\,dk.
$$
Thus, $K$ is proportional to $\om$.

\smallskip

{\it Step 2.} For $T_2$ we obtain
\be\la{T2}
T_2=-\int (e\we x)\rho\, dx= T_{21}+T_{22};\,\,\,\,%-\int (\hat e\we i\na_k)\hat\rho\, dk.
T_{21}:=-\int (e_{(0)}\we x)\rho\, dx,\,\,\,T_{22}:=-\int (e_{(r)}\we x)\rho\, dx.
\ee

Now we keep $T_{21}$ as it is and compute $T_{21}$ in Fourier space.

Recall that $\hat j(k,s)=\Om(s)\we i\na_k\hat\rho=i\Om(s)\we k\tilde\rho$. Then
$$
T_{22}(t):=-i\int dk\,\tilde\rho\left[k\we\int_0^t \frac{d}{dt}\hat K\vert_{t-s}\hat j(s)\,ds\right]=
-i\int_0^t\,ds\int dk\,\tilde\rho^2\cos|k|(t-s)[k\we(i\Om(s)\we k)]=
$$
$$
\int_0^t\,ds\int dk\,\tilde\rho^2\cos|k|(t-s)[\Om(s)k^2-k(k\cdot\Om(s))]=\int_0^t\,ds\,M_2(t-s)\Om(s),
$$
where
\begin{equation}\label{M2}
M_2(t-s):=\int dk\,\tilde\rho^2\cos|k|(t-s)[k^2E-k\otimes k].
\end{equation}

The matrix $k^2E-k\otimes k$ reads
$$
\left(
\ba{ccc}
k_2^2+k_3^2 & -k_1k_2 & -k_1k_3 \\
-k_2k_1 & k_1^2+k_3^2 & -k_2k_3 \\
-k_3k_1 & -k_3k_2 & k_1^2+k_2^2
\ea
\right).
$$
Thus, every non-diagonal element in (\ref{M2}) is zero, since the integrand function is odd w.r.t. corresponding variable. (Note that convergence of each integral above and below is provided by sufficient decay of $\hat\rho$ and its derivatives in $k$ due to (\ref{eq1}).) Further,
$$
\int dk\,\tilde\rho^2\cos|k|(t-s)k_j^2=\frac13\int dk\,\tilde\rho^2\cos|k|(t-s)k^2,\,\,\,j=1,2,3
$$
by the change of variables $k_i\mapsto k_j$ and we obtain
\begin{equation}\label{M2corr}
M_2(t-s):=\frac23\int dk\,\tilde\rho^2\cos|k|(t-s)k^2E.
\end{equation}
Finally, we change once more the order of integration and obtain
\begin{equation}\label{T22corr}
T_{22}(t)=\frac23\int dk\,\tilde\rho\,^2k^2\int_0^t\,ds\cos|k|(t-s)\Om(s).
\end{equation}

\smallskip

{\it Step 3.} $T_3=\int\,x\we[(\om(t)\we x)\we b]\rho\,dx=\om(t)\we\int\,(x\cdot b)\,x\rho\,dx=T_{31}+T_{32}$, where
\be\la{T3}
T_{31}:=\om(t)\we\int\,(x\cdot b_{(0)})\,x\rho\,dx,\,\,\,\,T_{32}:=\om(t)\we\int\,(x\cdot b_{(r)})\,x\rho\,dx.
\ee

We keep $T_{31}$ as it is and compute $T_{32}$ in Fourier space: $\,\,\,\,T_{32}(t):=$
$$
\om(t)\we\int\,dk\left[-\int_0^t\,ds\,im\hat K_{(t-s)}\hat j(k,s)\right]\tilde\rho=\om(t)\we\int_0^t\,ds\,\int\,dk\tilde\rho^2\left[\hat K_{(t-s)}k\we(\Om(s)\we k)\right]=
$$
\begin{equation}\label{T32corr}
\om(t)\we\int_0^t\,ds\,\int\,dk\,\tilde\rho^2\hat K_{(t-s)}(\Om(s)k^2-k(\Om(s)\cdot k))=\om(t)\we\frac23\int\,dk\,\tilde\rho^2k^2\int_0^t\,ds\,\frac{\sin|k|(t-s)}{|k|}\Om(s),
\end{equation}
similarly to the case of $T_{22}$.

Finally, we rewrite (\ref{eqOm2}) as
\begin{equation}\label{eqOm3f}
\dot\Om = M\Om + T(t),
\end{equation}
where $M$ is the skew-adjoint matrix of $m:=(K/I)\we$ and
\begin{equation}\label{TT}
T(t):=(1/I)(T_{21}(t)+T_{31}(t)+T_{22}(t)+T_{32}(t)).
\end{equation}

Note that, in detail,
$$
\frac{K}{I}=\frac{-\frac23{\om}\int\tilde\rho\Delta_k\hat\rho\,dk}{\frac23\int x^2\rho(x)\,dx}=\frac{-\frac23{\om}\int\tilde\rho\Delta_k\hat\rho\,dk}{-\frac23\int \Delta_k\hat\rho(k)\,dk} = \alpha(\rho)\om;\,\,\,\,
\alpha(\rho):=\frac{\int\tilde\rho\Delta_k\hat\rho\,dk}{\int \Delta_k\hat\rho(k)\,dk}.
$$

\begin{remark}
{\rm In \cite{I21}, by computational error it was claimed that $T_{22}=0$, $T_{23}=0$. The analysis below shows that this error does not affect the validity of the Theorem \ref{t1}.}
\end{remark}

\subsection{Integral inequality}

Let us write the integral equation equivalent to (\ref{eqOm3f}):

\be\la{integr}
\Om(t)=\Om_0+\alpha\om\we\int\limits_0^t\Om(\tau)\,d\tau+(1/I)\int\limits_0^t(T_{21}(\tau)+T_{31}(\tau)+T_{22}(\tau)+T_{32}(\tau))\,d\tau.
\ee

We estimate $|\Om(t)|$ by this equation:

\be\la{integr1}
|\alpha\om\we\int\limits_0^t\Om(\tau)\,d\tau|\le|\alpha\om|\cdot\int\limits_0^t|\Om(\tau)|\,d\tau.
\ee

\be\la{integr2}
|\int\limits_0^t(T_{21}(\tau)\,d\tau|\le C(\rho,e_0,b_0)\Vert(e_0,b_0)\Vert_{L^2\times L^2},
\ee
by (\re{T2}), (\re{zero}) and since $U(\tau)$ is isometric. Similarly,
\be\la{integr3}
|\int\limits_0^t(T_{31}(\tau)\,d\tau|\le C(\rho,e_0,b_0)\Vert(e_0,b_0)\Vert_{L^2\times L^2},
\ee
also by (\re{econ}). Now put $v(\tau):=\sup\limits_{[0,\tau]}|\Om(s)|$. Then
$$
|(1/I)\int\limits_0^t\,T_{22}(\tau)\,d\tau|\le\frac2{3I}\int\limits_0^t\int\,dkk^2\tilde\rho^2|\int\limits_0^\tau\,ds\cos(|k|(\tau-s))\Om(s)|\le
C(\rho)\int\limits_0^t\,\tau v(\tau)\,d\tau.
$$

Similarly,
$$
|(1/I)\int\limits_0^t\,T_{32}(\tau)\,d\tau|\le\frac2{3I}|\om(t)|\int\limits_0^t\,d\tau\tau v(\tau)\int\,dk\tilde\rho^2(k)\le  C(\rho,e_0,b_0)\int\limits_0^t\,d\tau\tau v(\tau).
$$

As the result, we get
\be\la{Omt}
|\Om(t)|\le||\Om_0|+C(\rho,e_0,b_0)\int\limits_0^t(\delta(e_0, b_0)+(1+\tau)v(\tau))\,d\tau.
\ee
Here $\delta(e_0, b_0)\to0$ as $\Vert(e_0,b_0)\Vert_{L^2\times L^2}\to0$.

Now we show that in the left hand side of (\re{Omt}) $|\Om(t)|$ can be replaced by $v(t)$. Indeed, by continuity, $v(t)=|\Om(t_0)|$, $t_0\in[0;t]$. Then
$$
v(t)=|\Om(t_0)|\le C(\rho,e_0,b_0)\int\limits_0^{t_0}(\delta(e_0, b_0)+(1+\tau)v(\tau))\,d\tau\le C(\rho,e_0,b_0)\int\limits_0^t(\delta(e_0, b_0)+(1+\tau)v(\tau))\,d\tau.
$$
By the strengthened Gronwall's lemma,
\be\la{stabzero}
v(t)\le C(\rho,e_0,b_0)\delta(e_0, b_0)t+|\Om_0|e^{C(\rho,e_0,b_0)(t+t^2/2)}.
\ee

From (\re{stabzero})
it follows that for any finite $T\in\R$
\begin{equation}\label{cdOm}
\Om(\cdot)\to0\,\,\,{\rm in}\,\,\, C(0,T;\R^3)\,\,\,{\rm as}\,\,\,\Om_0\to0\,\,\,{\rm in}\,\,\,\R^3,\,\,\,{\rm in\,\,\, particular,}\,\,\,\Om(T)\to0.
\end{equation}
By (\ref{cdOm}), (\ref{ek}), (\ref{bk})
\begin{equation}\label{cd}
(\Om(T),e(\cdot,T),b(\cdot,T))\to(0,0,0)\,\,\,{\rm in}\,\,\,L\,\,\,{\rm as}\,\,\,(\Om_0,e_0,b_0)\to(0,0,0)\,\,\,{\rm in}\,\,\,L.
\end{equation}

\subsection{The special case of zero initial fields perturbation}
Let us consider the special case of initial data $e_0(x)=0$, $b_0(x)=0$, and $\Om_0$ is arbitrary. In this case one has $e_{(0)}(x,t)=0$, $b_{(0)}(x,t)=0$ and hence, $T_{21}(t)=0$ and $T_{31}(t)=0$ by (\ref{T2}), (\ref{T3}). Then the equation (\ref{eqOm3f}) becomes a linear homogeneous integro-differential equation and its solution reads
\begin{equation}\label{linOm}
\Om(t)=A(t)\Om_0,
\end{equation}
where $A(t)$ is a $3\times3$-matrix and $\Vert A(t)\Vert$ is bounded uniformly in $t\in\R$ by (\ref{ub}).

Further, in this case
\begin{equation}\label{j}
j(x,s)=(\Om(s)\we x)\rho(x)=(A(s)\Om_0\we x)\rho(x).
\end{equation}
Then by %(\ref{Om0}) and
(\ref{U}), (\ref{j}), and (\ref{linOm}), %after the integration the right-hand side of (\ref{U0})
\begin{equation}\label{lineb}
\left(\begin{array}{c} e(x,t) \\ b(x,t) \end{array}\right) =  W(x,t)\Om_0\,,
\end{equation}
where $W(x,t)$ is a $6\times3$ matrix. The components $w_{ij}(t,x)$, $i=1, ... , 6$, $j=1, ... , 3$ of the matrix are functions bounded in $L^2$ uniformly in $t\ge0$ due to (\ref{ub}). %Remark 1
Then we obtain the following preliminary result on stability:

\begin{proposition}\label{p2}
For the system (\ref{M})-(\ref{lt}), in the phase space ${\cal M}$, the soliton \linebreak $(\om, E_\om, B_\om)$ is Lyapunov stable (and as well orbital stable) with respect to perturbations of type $(\Om_0, e_0=0, b_0=0)$.
\end{proposition}

\subsection{Completing the proof of Theorem \ref{t1}}

Consider initial perturbations $(e_0, b_0, \Om_0)$ such that
$$%\label{Om0}
{\rm supp}\,e_0\subset\{|x|\le R\},\,\,\,{\rm supp}\,b_0\subset\{|x|\le R\};
$$
with a fixed $R>0$.

By the strong Huygens principle for the group of the free Maxwell equations \cite{IKS11,KS00,IKS02,KS06} the supports of $e_{(0)}(x,t)$ and of $b_{(0)}(x,t)$ are subsets of the region $\{|x|>t-R\}$.

Then, since $\rho$ is compact supported, there is a $\overline T=\overline T(R,R_\rho)$ such that $T_{21}(t)=0$ and $T_{31}(t)=0$ for $t\ge\overline T$ by (\ref{zero}), (\ref{T2}), (\ref{T3}). Then for $t\ge\overline T$ the equation (\ref{eqOm3f}) reads
\begin{equation}\label{OmOm}
\dot\Om=M\Om+\frac1I(T_{22}(t)+T_{32}(t))\,\,\,{\rm with\,\, the\,\, initial\,\, condition}\,\,\, \overline{\Om}:=\Om(\overline T),
\end{equation}
where $\Om(t)$ is the solution to (\ref{eqOm3f}) for $0\le t\le\overline T$. The equation (\ref{OmOm}) is a linear homogeneous integro-differential equation w.r.t. $\Om$.

By (\ref{cdOm})
\begin{equation}\label{OmT}
\overline\Om\to0\,\,\,{\rm in}\,\,\,\R^3\,\,\,{\rm as}\,\,\,(\Om_0,e_0, b_0)\to0\,\,\,{\rm in}\,\,\,L.
\end{equation}

Further, for $t\ge\overline T$ the solution in (\ref{U}) reads
\begin{equation}\label{Ult}
\left(\ba{c} e(x,t) \\ b(x,t) \ea\right) = U(t)\left(\ba{c} e_0(x) \\ b_0(x) \ea\right)-
\int_0^{\overline T}U(t-s)\left(\ba{c} j(x,s) \\ 0 \ea\right)\,ds-\int_{\overline T}^tU(t-s)\left(\ba{c} j(x,s) \\ 0 \ea\right)\,ds.
\end{equation}

i) In the right-hand side of (\ref{Ult}), for the first term $(e_{(0)},b_{(0)})$ we have
\begin{equation}\label{unit}
\Vert e_{(0)}(\cdot,t)\Vert_{L^2}^2+\Vert b_{(0)}(\cdot,t)\Vert_{L^2}^2 = \Vert e_0\Vert_{L^2}^2+\Vert b_0\Vert_{L^2}^2,
\end{equation}
since the group $U(t)$ is unitary.

ii) For the second term let us observe that $s\in[0;\overline{T}]$ with a fixed $\overline{T}$, the group $U(t-s)$ is unitary, and $j(x,s)=(\Omega(s)\wedge x)\rho(x)$. Then, by continuous dependence, see  (\ref{cdOm}), %(\ref{cd}),
this term is bounded in $[L^2(\R^3;\R^3)]^2$ uniformly in $t\in\R$ and tends to zero as $(\Om_0,e_0, b_0)\to0$.

iii) For the third term it follows from i), ii), and (\ref{ub}) that it is also bounded in $[L^2(\R^3;\R^3)]^2$ uniformly in $t\in\R$. Further, the current reads
$$%\label{Om0}
j(x,s)=(A(s)\overline\Om\wedge x)\rho,\,\,\,s\ge\overline{T},
$$
and after integrating the term becomes $\overline{W(x,t)}\,\overline\Om$, where the components

\begin{equation}\label{third}
\overline{w}_{ij}(x,t)\,\,\,{\rm of\,\,the\,\,matrix}\,\,\,\overline{W}(x,t)\,\,\,{\rm are\,\,bounded\,\,in}\,\,\,L^2\,\,\,{\rm uniformly\,\,in}\,\,\,t,
\end{equation}
because of the uniform boundness (\ref{ub}). By (\ref{OmT}) this term tends to zero in $L$ as $(\Om_0,e_0,b_0)\to(0,0,0)$ in $L$.

iv) Finally, for $\Om(t)$ itself we have $\Om(t)\to0$ in $\R^3$ as $(\Om_0,e_0,b_0)\to(0,0,0)$ in $L$ by (\ref{cdOm}) for $t\le\overline{T}$ and by (\ref{linOm}) with $\overline\Om$ instead of $\Om_0$ for $t\ge\overline{T}$.

Then the conclusion of Theorem \ref{t1} follows from i) to iv). The proof is complete.

%%%%%%%%%%%%%%%%%%%%%%%%%%%%%%%%%%%%%%%%%%%%%%%%%%%%%
%%%%%%%%%%%%%%%%%%%%%%%%%%%%%%%%%%%%%%%%%%%%%%%%%%%%%

\section{Absence of attraction to a soliton of finite angular\\ momentum}

\setcounter{equation}{0}

In this section we show that there is no attraction, in the energy norm (\re{Lnorm}), to a soliton of finite angular momentum, for some solutions with initial data on the surface of states of the same angular momentum.

\subsection{Angular momentum}

The angular momentum is defined by
\be\la{am}
M(\om,E,B):=I\om+\int\,x\we(E(x)\we B(x))\,dx.
\ee

Note that for $(E,B)\in (L^2(\R^3;\R^3), L^2(\R^3;\R^3))$ the angular momentum is generally not defined. It is well defined for the fields $(E,B)$ with the finite weighted norms
\be\la{weight}
\int\,|x|\,|E(x,t)|^2dx,\,\,\,\,\int\,|x|\,|B(x,t)|^2dx.
\ee

\subsection{Faster spatial decay and angular momentum of solitons}

From (\re{solf}), by a straightforward computation we obtain that
$$
\pa_j\hat E_\om\in L^2,\,\,\,\,\pa_j\hat B_\om\in L^2,\,\,\,\,j=1,2,3
$$
under the condition
\be\la{rhozk}
\hat\rho(0)=0.
\ee
Hence,
\be\la{solweight}
x_j E_\om\in L^2,\,\,\,\,x_j B_\om\in L^2,\,\,\,\,j=1,2,3\,\,\,\,{\rm (in\,\,\,}x-{\rm space)}
\ee
under the same condition (\ref{rhozk}) or the equivalent condition
\be\la{rhozx}
\int\,\rho(x)\,dx=0.
\ee
As the result, we obtain the following

\begin{proposition}\label{p3}
Under the condition (\ref{rhozk}) or the equivalent condition (\ref{rhozx}) the weighted norms of the soliton fields
\be\la{weightfin}
\int\,|x|\,|E_\om(x)|^2dx,\,\,\,\,\int\,|x|\,|B_\om(x)|^2dx
\ee
are finite and the angular momentum of the soliton
\be\la{amts}
M_\om:=I\om+\int\,x\we(E_\om(x)\we B_\om(x))\,dx<\infty
\ee
is well-defined.
\end{proposition}

\begin{remark}
The angular momentum of a soliton is computed exactly in Appendix.
\end{remark}

%%%%%%%%%%%%%%%%%%%%%%%%%%%%%%%%%%%%%%%%%%%%%%%%%%%%%%%%%%%%%%%%%%%%%%%%%%%%%%%%%%%%%%%

\subsection{Partial negative result on attraction}

\begin{proposition}\label{t2}
Let the conditions (\re{rhozk}) hold. Consider the soliton $S_\om=(\om,E_\om,B_\om)$ of the finite angular momentum $M=M_\om$. Then there exists an initial state $(\Om,E,B)$ of the angular momentum $M$ such that for the solution with the initial condition $(\Om,E,B)$ there is no attraction to $S_\om$.
%Then for any $\ve>0$ the set ${\cal S}$ is not $\ve$-attractive.
\end{proposition}

Let us note that by the energy conservation a solution $Y(t)$ with initial data $Y_0=(\om_0,E_0,B_0)$ cannot tend to a soliton $S_\om=(\om, E_\om,B_\om)$ if $H(S_\om)\ne H(Y_0)$.

On the other hand,
\be\la{HSom}
H(S_\om)=I\om^2/2+(1/2)\int\,(|E_\om|^2+|B_\om|^2)dx.
\ee
For $\om=0$ one has
$$
H_0:=H(S_0)=(1/2)\int\,|E_0|^2dx=(1/2)\int\,|\hat E_0|^2dk =(1/2)\int\,(\hat\rho^2/k^2)dk
$$
by (\re{solf}). This is the minimal value of $H(S_\om)$ which is reached at the zero soliton.

If one could construct initial data with an energy less than $H_0$, this would mean that there is no global attraction to the set of all solitons. But in fact the minimal possible value for initial data is exactly $H_0$. To prove this one has to set the conditional extremum problem $\Vert E\Vert^2_{L^2}\to$~min, $\na\cdot E=\rho$. By the standard Lagrange method we obtain that the minimal value is $H_0$.

Further, also by (\re{solf}), (\re{HSom}) $H(S_\om)$ is continuous in $\om$ and takes any value of $[H_0;+\infty)$. Thus, energy argument does not interfere the global attraction. That is why we apply the angular momentum argument to obtain the partial negative result on attraction.

%%%%%%%%%%%%%%%%%%%%%%%%%%%%%%%%%%%

\subsection{Proof of Proposition \ref{t2}}

\subsubsection{Energy variation on the surface of constant angular momentum}

Let us fix an arbitrary $\om$ and consider the soliton $(\om,E_\om,B_\om)$ with the angular momentum $M=M_\om$. Consider the surface of constant angular momentum
$$
{\cal S}_\om:=\{(\Om,E,B): M(\Om,E,B)=M=M_\om\}.
$$
 where the fields $(E,B)$ obey the constraints (\re{div}). Let us make variations of the energy $H(\Om,E,B,)$ on this surface. For this purpose we express $\Om$ in $M$,
\be\la{Omega}
\Om=\fr1I(M-\int x\we(E\we B)dx)
\ee
and make variations of
$$
H_M(E,B)=\fr1{2I}\left(M-\int x\we(E\we B)dx\right)^2+\fr12\int(|E(x)|^2+|B(x)|^2)dx
$$
in $E$ and $B$.

Let $e$, $b$ be variations in $E,B$ of finite weighted norms, such that $H_M(E+e,B+b)$ is finite. Since $\na\cdot E=\rho$ and $\na\cdot(E+e)=\rho$, we have \be\la{dive}\na\cdot e=0.\ee Similarly \be\la{divb}\na\cdot b=0.\ee
Let us find variation in $E$,
$$
\fr{d}{dt}\Big|_{t=0}H_M(E+te,B).
$$
We have
$$
\fr{d}{dt}\Big|_{t=0}\left[\fr1{2I}\left(M-\int x\we((E+te)\we B)dx\right)^2+\fr12\int(|E+te|^2+|B|^2)dx\right]=
$$
$$
\int(E+te)\cdot e\,dx\Big|_{t=0}-\fr1{I}\left(M-\int x\we((E+te)\we B)dx\right)\Big|_{t=0}\cdot\int x\we(e\we B)dx=
$$
\be\la{varE}
\int E\cdot e\,dx-\fr1{I}\left(M-\int x\we(E\we B)dx\right)\cdot\int x\we(e\we B)dx.
\ee
Similarly, variation in $b$ results in
\be\la{varB}
\int B\cdot b\,dx-\fr1{I}\left(M-\int x\we(E\we B)dx\right)\cdot\int x\we(E\we b)dx.
\ee

%%%%%%%%%%%%%%%%%%%%%%%%%%%%%%%%%%%%%%%%%%%%%%%%%%%%%%%%%%%%%%%%%%%%%%%%%%%%%%%%%%%%%%%

\subsubsection{Absence of attraction to the soliton on the surface of constant angular momentum. Energy increment %alternation
argument}

For the absence of attraction to $S_\om$ it is sufficient to show that the energy $H_M$ is not constant on the surface ${\cal S}_\om$ of constant angular momentum. Indeed, if at some point of ${\cal S}_\om$ the energy differs from that of the soliton, the solution starting at that point cannot tend to the soliton by energy conservation.

Let us start with an arbitrary point $(\Om,E,B)$, $B\ne0$ of ${\cal S}_\om$ and put $e=0$, $b=B$. Then the variation (\re{varE}) in $E$ at the point $(\Om,E,B)$ equals zero, while the variation (\re{varB}) in $B$ equals
$$
V_B(\Om,E,B)=\int\,dx\,\left(B^2+\fr1I(M_0^2-MM_0)\right),\,\,\,{\rm where}\,\,\,M=M_\om,\,\,\,M_0=\int\,dx\,x\we(E\we B).
$$
Note that the first term $\int\,dx\,B^2$ is positive and does not depend on $E$. If $M_0=0$, $V_B(E,B)>0$. If $M_0\ne0$, we can make $V_B(E,B)>0$ as well changing $E$ to $\lambda E$ with sufficiently large $\lambda$. Thus, at a certain point of ${\cal S}_\om$ the energy $H_M$ locally increases in the direction of the tangent vector $b=B$ to ${\cal S}_\om$ and hence is not constant.

Proposition is proved.

%%%%%%%%%%%%%%%%%%%%%%%%%%%%%%%%%%%%%%%%%%%%%%%%%%%%%%%%%%%%%%%%%%%%%%%%

\section{The absence of attraction in 2D case}

\setcounter{equation}{0}

The 2D Maxwell-Lorentz system corresponding to (\re{odd-even})--(\re{lt}) reads, \ci{KKrest2D}:

\be\la{ml2d}
\left\{\ba{l}
\dot E(x;t) = J\na B(x;t) + \om(t)Jx\rho(x),\\
\dot B(x;t) = -\na\cdot JE(x;t),\,\,\na E(x;t) = \rho(x),\\
I\dot\om(t) = \int dx x\cdot JE(x;t)\rho(x).
\ea\right.
\ee
Here $x=(x_1,x_2)\in\R^2$, $t\in\R$, $E(x,t)=(E_1(x,t),E_2(x,t))$ is $\R^2$-valued, $B(x,t)$ and $\om$ are $\R$-valued;
$$
J=\left(
\ba{ll}
0 & 1\\
-1 & 0
\ea
\right)
$$

The existence of dynamics for the system (\re{ml2d}) is established in \ci{KKrest2D}. The energy
\be\la{E2D}
H(\om,E.B)=\fr{I}{2}\om^2+\fr12\int(|E(x)|^2+|B(x)|^2)dx
\ee
is conserved along the solutions. The solitons, in Fourier space read
\be\la{sol2d}
\hat E_\om=-\fr{i\hat\rho k}{k^2},\,\,\,\,\hat B_\om=-\fr{\om k\cdot\na\hat\rho}{k^2}.
\ee

Let $\rho$ satisfy the conditions
\begin{equation}\label{eq12d}
\rho\in C_0^{\infty}(\R^2),\,\,\,\,\rho(x)= \rho_{rad}(|x|),\,\,\,\,\rho(x)=0\,\,\,{\rm for}\,\,\,|x|>R_{\rho}>0.\,\,\,\hat\rho(0)=0. %\eqno{(C)}
\end{equation}
Then the angular momentum of solitons is well defined end equals, similarly to (\re{amsfin}),
\be\la{amsfin2d}
M_\om:=I\om+\int x\we(E_\om\we B_\om)dx=I\om-\int (x\cdot E_\om)B_\om\,dx=\om\left(I+\int\fr{|\na\hat\rho|^2}{r^2}dk\right)
\ee
which coincides with \ci[(2.14)]{KKrest2D}. Let us fix an arbitrary $\om$ and consider the soliton $(\om,E_\om,B_\om)$ with the angular momentum $M=M_\om$. On the surface of constant angular momentum
$$
{\cal M}_\om:=\{(\Om,E,B): M(\Om,E,B)=M=M_\om\}
$$
with $\na\cdot E=\rho$ we make, in $E$ and $B$, variations of the energy $H(\Om,E,B,)=H_M(E,B)$, where
$$
H_M(E,B)=\fr1{2I}\left(M+\int (x\cdot E)B\,dx\right)^2+\fr12\int(|E(x)|^2+|B(x)|^2)dx.
$$
The variations in $E$ and $B$ equal correspondingly ($\na\cdot e=0$, $\na\cdot b=0$)
$$
V_E(\Om,E,B)=\int E\cdot e\,dx+\fr1{I}\left(M+\int (x\cdot E)B\,dx\right)\cdot\int (x\cdot e)B\,dx,
$$
$$
V_B(\Om,E,B)=\int B\cdot b\,dx+\fr1{I}\left(M+\int (x\cdot E)B\,dx\right)\cdot\int (x\cdot E)b\,dx.
$$
Similarly to Subection 4.4.2 we find a point in ${\cal M}_\om$ such that the energy $H_M$ locally increases at this point and hence is not constant on ${\cal M}_\om$. Then there exist initial data $(\Om, E, B)\in {\cal M}_\om$ such that the solution with these initial data does not tend to $S_\om$ in the energy norm corresponding to (\re{E2D}).

We obtain the following statement for the 2D case:

\begin{proposition}\label{t22D}
Let for the 2D system (\re{ml2d}) the conditions (\re{eq12d}) hold. Consider the soliton $S_\om=(\om,E_\om,B_\om)$ of the finite angular momentum $M=M_\om$. Then there exists an initial state $(\Om,E,B)$ of the angular momentum $M$ such that for the solution with the initial condition $(\Om,E,B)$ there is no attraction, in the energy norm, to $S_\om$.
\end{proposition}

\subsection{Remark on global attraction for potential form of 2D Maxwel-Lorentz equations}

In \ci{KKrest2D} the 2D Maxwell-Lorentz system (\re{ml2d}) is transformed to potential form via

\be\la{ml2dptrans}
E=-\dot A-\na\Phi,\,\,\,B=\na\cdot(JA),\,\,\,\na\cdot A=0,\,\,\,\De\Phi=-\rho.
\ee
Here $A=(A_1(x_1,x_2,t), A_2(x_1,x_2,t)$ and $\Phi(x_1,x_2)$ is determined by the last equation of (\re{ml2dptrans}).

The system (\re{ml2d}) in potential form reads

\be\la{ml2dpot}
\left\{
\ba{l}
\dot A=\Pi\\
\dot\Pi=\De A-\om Jx\rho\\
I\dot\om=-\int dx\, J\Pi\cdot x\rho
\ea
\right.
\ee

The energy

\be\la{epot}
H_{\rm p}=\fr{I\om^2}{2}+\fr12\int dx (|\na A|^2+|\Pi|^2)
\ee
and the angular momentum
\be\la{Mpot}
M_{\rm p}=I\om-\int dx (A\cdot Jx\rho)
\ee
are conserved along solutions to (\re{ml2dpot}). A soliton, in terms of potentials reads
$$
S_\om=(A_\om(x),0),\,\,\,\,\hat A_\om(k)=-\fr{\om J\widehat{x\rho}(k)}{k^2}.
$$

Denote
\be\la{Mfield}
M_{\rm f}:=I\om+\int dx\,x\we(E\we B)
\ee
the angular momentum expressed in terms of fields. At a soliton $S_\om$, $M_{\rm f}=M_{\rm p}$, but at an arbitrary point $(\om,A,\Pi)$ the values of $M_{\rm f}$ and $M_{\rm p}$ can differ. Indeed, $M_{\rm f}$ expressed in terms of potentials reads
\be\la{Mfp}
M_{\rm f}=I\om-\int dx\,x\cdot(\Pi+\na\Phi)(\na\cdot JA)
\ee
which depends on $\Pi$ while $M_{\rm p}$ does not depend on $\Pi$. Then varying $\Pi$ we can obtain different values of $M_{\rm p}$ and $M_{\rm f}$ for the same point $(\om,A,\Pi)$.

Anyway, consider the energy $H_{\rm p}$ on the surface of constant angular momentum $M_{\rm p}$, the energy reads
\be\la{epMp}
H_{{\rm p}, M}=\fr12\int dx (|\na A|^2+|\Pi|^2)+\fr1{2I}\left(M+\int dx (A\cdot Jx\rho)\right)^2.
\ee

Variation of $H_{{\rm p}, M}$ in $A$ equals
\be\la{VAepMp}
V_A(H_{{\rm p}, M}):=\fr{d}{dt}\vert_{t=0}H_{{\rm p}, M}(A+ta,\Pi)=\int dx \na A\cdot\na a+\fr1{I}(M+\int dx (A\cdot Jx\rho))\cdot\int dx (a\cdot Jx\rho).
\ee
Variation in $\Pi$ equals
\be\la{VPiepMp}
V_\Pi(H_{{\rm p}, M}):=\fr{d}{dt}\vert_{t=0}H_{{\rm p}, M}(A,\Pi+t\pi)=\int dx \Pi\cdot\pi.
\ee

For $a=0$, $\pi=\Pi$ we have $V_A(H_{{\rm p}, M})=0$, $V_\Pi(H_{{\rm p}, M})=\int dx |\Pi|^2>0$ if $\Pi\ne0$. Thus, there is a point, where $H_{{\rm p}, M}$ locally increases in the direction $\Pi$ and hence is not constant on the surface of constant angular momentum. We obtain the ``potential'' analogue of Proposition 5.1:

\begin{proposition}\label{t22Dp}
Let for the 2D system (\re{ml2dpot}) the conditions (\re{eq12d}) hold. Consider the soliton $S_\om=(\om,A_\om,0)$ of the finite angular momentum $M=M_\om$. Then there exists an initial state $(\Om,A,\Pi)$ of the angular momentum $M$ such that for the solution with the initial condition $(\Om,A,\Pi)$ there is no attraction, in the energy norm corresponding to (\re{epMp}), to $S_\om$.
\end{proposition}

It is interesting to note that in \ci{KKrest2D} it is shown that for initial data of finite angular momentum the solution to the system (\re{ml2dpot}) converges to the soliton of the same angular momentum, \ci[Theorem 3.2]{KKrest2D}. This result does not contradict to Proposition 5.2 because in \ci{KKrest2D} the convergence is established in weak weighted norm \ci[Definition 3.1]{KKrest2D} with the weight parameter $\beta<-5/2$. While our argument shows that the convergence in the energy norm does not take place.

\medskip

{\bf Gratitudes.} The author is grateful to Professor Herbert Spohn, Technical University of Munich, Germany, for drawing the author's attention to the problem, Professor Alexander Komech, Moscow University, Russia and University of Vienna, Austria, and Dr. Elena Kopylova, University of Vienna, Austria, for fruitful discussions.

\section{Appendix. Angular momentum of a soliton}

\setcounter{equation}{0}

Let us compute the angular momentum $M_\om$ of a soliton $(E_\om,B_\om,\om)$ given by (\re{amts}), in detail. We have
$$
\int x\we(E_\om\we B_\om)dx=\int(E_\om(x\cdot B_\om)-B_\om(x\cdot E_\om))dx=
$$
\be\la{amsolit}
\int\left[(i\na_k\cdot\hat B_\om)\ov{\hat E_\om}-\hat B_\om\ov{(i\na_k\cdot\hat E_\om)}\right]dk=
i\int\left[(\na_k\cdot\hat B_\om)\ov{\hat E_\om}+\hat B_\om(\na_k\cdot\ov{\hat E_\om})\right]dk.
\ee
Since
$$
\hat B_\om=\fr{\rho_r'((k\om)k-k^2\om)}{r^3},\,\,\,\,\ov{\hat E_\om}=\fr{ik\ov{\ti\rho}}{k^2},
$$
we obtain by the computation of Section 4.1 that
$$
\na_k\cdot\hat B_\om=\fr{2\rho_r'(k\om)}{r^3},\,\,\,\,\na_k\cdot\ov{\hat E_\om}= i\fr{\rho_r+\rho_r'r}{r^2}.
$$
Then the expression in (\re{amsolit}) equals
$$
-\int\left[\fr{2\rho_r'\rho_r(k\om)k}{r^5}+((k\om)k-k^2\om)\rho_r'\fr{\rho_r+\rho_r'r}{r^5}\right]dk.
$$
Since $\rho_r$ and $\rho_r'$ are even in $k$, the last expression simplifies to
$$
\fr23\om\int\fr{(\rho'_r)^2}{r^2}dk=\fr23\om\int\fr{|\na\hat\rho|^2}{r^2}dk.
$$
Finally, the angular momentum of the soliton equals
\be\la{amsfin}
M_\om=\om\left(I+\fr23\int\fr{|\na\hat\rho|^2}{r^2}dk\right).
\ee

\br In $x$-representation (\re{amsfin}) reads
$$
M_\om=\om\left(I+\fr23\int\,x\rho\cdot\De^{-1}(x\rho)dx\right).
$$
\er

\end{document}